\documentclass[11pt]{article}
\usepackage{amsfonts,epsf,amsmath,amssymb}
\usepackage{tikz, caption, subcaption}
\usepackage[numbers,sort&compress]{natbib}

\newtheorem{theorem}{\bf Theorem}[section]
\newtheorem{corollary}[theorem]{\bf Corollary}

\newtheorem{lemma}[theorem]{\bf Lemma}
\newtheorem{proposition}[theorem]{\bf Proposition}
\newtheorem{conjecture}[theorem]{\bf Conjecture}

\newcommand{\vertex}{\node[vertex]}
\tikzstyle{vertex}=[circle, draw, inner sep=0pt, minimum size=6pt]
\newcommand{\pch}{\chi_{\rho}}
\newcommand{\diam}{{\rm diam}}
\newcommand{\proof}{\noindent{\bf Proof.\ }}
\newcommand{\qed}{\hfill $\square$ \bigskip}
\newcommand{\smallqed}{{\tiny ($\Box$)}}

\begin{document}

\title{\bf Packing chromatic number, $(1,1,2,2)$-colorings, and characterizing the Petersen graph}

\author{
 \phantom{xxxxx} \and Bo\v{s}tjan Bre\v{s}ar $^{a,b}$  \and Sandi Klav\v zar $^{a,b,c}$ \and \phantom{xxx}
\and Douglas F. Rall $^{d}$ \and Kirsti Wash $^e$}

\date{}

\maketitle

\begin{center}
$^a$ Faculty of Natural Sciences and Mathematics, University of Maribor, Slovenia\\
\medskip

$^b$ Institute of Mathematics, Physics and Mechanics, Ljubljana, Slovenia\\
\medskip

$^c$ Faculty of Mathematics and Physics, University of Ljubljana, Slovenia\\
\medskip

$^d$ Department of Mathematics, Furman University, Greenville, SC, USA\\
\medskip

$^e$ Department of Mathematics, Trinity College, Hartford, CT, USA\\
\end{center}

\begin{abstract}
The packing chromatic number $\pch(G)$ of a graph $G$ is the smallest integer $k$ such that the vertex set of $G$ can be partitioned into sets $\Pi_1,\ldots,\Pi_k$, where $\Pi_i$,  $i\in [k]$, is an $i$-packing. The following conjecture is posed and studied: if $G$ is a subcubic graph, then $\pch(S(G))\le 5$, where $S(G)$ is the subdivision of $G$. The conjecture is proved for all generalized prisms of cycles. To get this result it is proved that if $G$ is a generalized prism of a cycle, then $G$ is $(1,1,2,2)$-colorable if and only if $G$ is not the Petersen graph. The validity of the conjecture is further proved for graphs that can be obtained from generalized prisms in such a way that one of the two $n$-cycles in the edge set of a generalized prism is replaced by a union of cycles among which at most one is a 5-cycle. The packing chromatic number of graphs obtained by subdividing each of its edges a fixed number of times is also considered.
\end{abstract}


\medskip\noindent
{\bf AMS Subj.\ Class:} 05C70, 05C15, 05C12


\section{Introduction}
Given a graph $G$ and a positive integer $i$, an {\em $i$-packing} in $G$ is a subset $W$ of the vertex set of $G$ such that the distance between any two distinct vertices from $W$ is greater than $i$. This generalizes the notion of an independent set, which is equivalent to a $1$-packing. Now, the {\em packing chromatic number} of $G$ is the smallest integer $k$ such that the vertex set of $G$ can be partitioned into sets $\Pi_1,\ldots,\Pi_k$, where $\Pi_i$ is an $i$-packing for each $i\in [k]$. This invariant is well defined on any graph $G$ and is denoted $\pch(G)$. It was introduced in~\cite{goddard-2008} under the name broadcast chromatic number, and subsequently studied under the current name, see~\cite{argiroffo-2014, bkr-2007, bkr-2016, bkrw-2016+,  ekstein-2014, fiala-2010, fiala-2009, finbow-2010, jacobs-2013, korze-2014, togni-2014, torres-2015}.

One of the intriguing problems related to the packing chromatic number is whether it is bounded by a constant in the class of all cubic graphs. In particular, it was asked already in the seminal paper~\cite{goddard-2008} what is the maximum packing chromatic number in the class of cubic graphs of a given order. Gastineau and Togni found a cubic graph $G$ with $\pch(G)=13$ and asked whether $13$ is an upper bound for $\pch$ in the class of cubic graphs~\cite{gt-2016}, which we answered recently in the negative~\cite{bkrw-2016+}.
More specifically, it was asked in~\cite{fiala-2009} whether the invariant is bounded in the class of  planar cubic graphs. A question of similar nature from~\cite{gt-2016} asks whether the subdivision $S(G)$ of any subcubic graph $G$ (i.e., a graph with maximum degree 3) has packing chromatic number no more than 5. This question is the main motivation for the present paper. We suspect that the answer is positive, and pose it as the following conjecture.

\begin{conjecture}
\label{con:main}
If $G$ is a subcubic graph, then $\pch(S(G))\le 5$.
\end{conjecture}

The packing chromatic number of subdivided graphs has been studied in several papers. Using subdivided graphs the class of graphs with  packing chromatic number equal to $3$ was characterized in~\cite{goddard-2008}. The effect on the invariant of the subdivision of an edge of a graph was analyzed in~\cite{bkrw-2016+}. It was observed in~\cite{bkr-2007} that $\pch(S(G))\le \pch(G)+1$ for any graph $G$, and further proved that $\pch(S(K_n))=n+1$. Consequently the packing chromatic number of subdivided graphs is generally not bounded, hence the restriction to subcubic graphs in Conjecture~\ref{con:main} is natural.

The paper is organized as follows. In the next section we introduce notation needed, list several facts related to Conjecture~\ref{con:main}, and prove a connection between the packing chromatic number (of subdivided graphs) and the so-called $(1,1,2,2)$-colorings. This connection is then used as our main tool while attacking the conjecture. Then, in Section~\ref{sec:generalized-prisms}, we prove that Conjecture~\ref{con:main} holds true for all generalized prisms of cycles. Along the way a characterization of the Petersen graph is obtained. (We refer to~\cite{goddard-2009} for a recent characterization of the Petersen graph and to~\cite{torri-1997} for older characterizations.) Moreover, it is shown that any optimal packing coloring of the subdivided Petersen graph looks differently than one would expect. In Section~\ref{sec:more-general} we then extend the main result of the previous section to the graphs obtained from generalized prisms in such a way that one of the two $n$-cycles in the edge set of a generalized prism is replaced by a union of cycles among which at most one is a 5-cycle. In the final section we consider the packing chromatic number of graphs obtained by subdividing each of its edges a fixed number of times.

\section{Notation and preliminary results}
\label{sec:preliminaries}

All graphs considered in this paper are simple and connected, unless stated otherwise.

Let $G$ be a graph and $S(G)$ its {\em subdivision}, that is, the graph obtained from $G$ by replacing each edge with a disjoint
 path of length 2. In other words, $S(G)$ is obtained from $G$ by subdividing each edge $e$ of $G$ with a new vertex to be denoted by $v_e$.
 The resulting vertex set $V(S(G))$ can thus  be considered as $V(G)\cup \{v_e \mid e\in E(G)\}$. More generally, if $i\ge 1$, we
 define the graph $S_i(G)$ as the graph obtained from $G$ by subdividing each of its edges precisely $i$ times. In other words, $S_i(G)$
 is obtained from $G$ by replacing each edge with a disjoint path of length $i+1$. Note that $S_1(G)=S(G)$.

Observe that if $H$ is a subgraph of $G$, then $\pch(H)\le \pch(G)$. Indeed, this follows because $d_H(u,v)\ge d_G(u,v)$ holds for any
vertices $u,v\in V(H)$. Consequently, a packing coloring of $G$ restricted to $H$ is  a packing coloring of $H$. Since every subcubic
graph is a subgraph of a cubic graph (easy exercise), it suffices to prove Conjecture~\ref{con:main} for cubic graphs. In addition, the following fact is a consequence of the characterization of the graphs of packing chromatic number 3 from~\cite{goddard-2008}.

\begin{proposition}
[\cite{goddard-2008}]
\label{prp:bipartite}
If $G$ is a (connected) bipartite graph of order at least $3$, then $\pch(S(G))=3$.
\end{proposition}

Hence we can restrict our attention to cubic non-bipartite graphs. Since $\pch(S(K_4))=5$ (see~\cite{bkr-2007}), Conjecture~\ref{con:main} reduces to 3-chromatic cubic graphs.  Before we continue, we demonstrate that the conjecture does not hold for all 3-chromatic graphs.

\begin{proposition}
If $K_{n,n,n}$ is the complete tripartite graph with all parts of order $n$,
then $\pch(S(K_{n,n,n}))\xrightarrow[n\to \infty]{}\infty$.
\end{proposition}

\proof
Let $G_n$ denote $S(K_{n,n,n})$. Since $\diam(G_{n})=4$ for $n\ge 2$, we infer that in any packing coloring $c$ of $G_n$ every color bigger
than $3$ appears at most once. Let $A,B$ and $C$ be the tripartition of $V(K_{n,n,n})$. Suppose there is a vertex $x$ from $A\cup B\cup C$
with $c(x)=1$. Since $N[x]$ induces $K_{1,2n}$ with $x$ as its center, in this case $c$ uses at least $2n$ colors. Otherwise, we may assume
without loss of generality that in $A\cup B$ there are vertices $y$ and $z$ with $c(y)=2$, $c(z)=3$. Clearly, then no vertex from $C$
can receive colors $2$ or $3$, which in turn implies that $c$ uses $n$ different colors on $C$.
\qed

Note that $\pch$ can be defined also in terms of a function on the vertex set of a graph $G$. Indeed, we say that a function $c:V(G) \to [k]$ is a {\em $k$-packing coloring} of $G$ if for each $i$ from the range of $c$, the set $c^{-1}(i)$ is an $i$-packing in $G$; we then also say that $G$ is {\em $k$-packing colorable}. In this way, $\pch(G)$ is the smallest integer $k$ such that there exists a $k$-packing coloring of $G$.

One approach to attack Conjecture~\ref{con:main} is by using the concept of an $S$-coloring, which generalizes that of a packing coloring. This concept was first briefly mentioned in~\cite{goddard-2008} and later formally introduced in~\cite{goddard-2012} as follows. Given a graph $G$ and a non-decreasing sequence $S=(s_1,\ldots,s_k)$ of positive integers, an {\em $S$-coloring} of $G$ is a partition of the vertex set of $G$ into $k$ subsets $\Pi_1,\ldots,\Pi_k$, where $\Pi_i$ is an $s_i$-packing for each $i\in [k]$. We say that $G$ is {\em $S$-colorable} if it has an $S$-coloring. Clearly, $\pch(G)\le k$ if and only if $G$ is $S$-colorable for $S=(1,2,\ldots,k)$. For further results on the $S$-packing coloring see~\cite{gastineau-2015a, gastineau-2015b, goddard-2014}.

The following result shows in what way $(1,1,2,2)$-colorable graphs are related to Conjecture~\ref{con:main}.

\begin{proposition} \label{prop:necessary}
If $G$ is $(1,1,2,2)$-colorable, then $\pch(S(G))\le 5$.
\end{proposition}

\proof
By \cite[Proposition 1]{gt-2016}, every $(1,1,2,2)$-colorable graph $G$ yields a $(1,3,3,5,5)$-colorable $S(G)$,
which in turn implies that $S(G)$ is $(1,2,3,4,5)$-colorable, that is, $S(G)$ is $5$-packing colorable.
\qed

We next state a result that will be the main tool in our subsequent proofs. For its statement recall that the {\em square} $G^2$ of a graph $G$ is the graph having the same vertex set as $G$ and two vertices are adjacent in $G^2$ precisely when their distance in $G$ is at most $2$.

\begin{lemma}\label{lem:CP1}
A graph $G$ is $(1,1,2,2)$-colorable if and only if there is a partition $\{V_1, V_2, V_3\}$ of $V(G)$
such that $V_2$ and $V_3$ are independent sets and $V_1$ induces a bipartite graph in $G^2$.
\end{lemma}
\proof
Suppose that $\{V_1, V_2, V_3\}$ is a partition of $V(G)$ as stated above. Let $A$ and $B$ represent the partite sets of the graph $G^2[V_1]$. Note that
$A$ is a $2$-packing in $G$ for otherwise $A$ would not be an independent set in $G^2$. Similarly, $B$ is a $2$-packing. Construct a $(1,1,2,2)$-coloring
of $G$ by assigning all the vertices of $V_2$ color $1$, all the vertices of $V_3$ color $2$, all the vertices of $A$ color $3$ and all the vertices of $B$ color $4$.
Thus, $(V_2,V_3,A,B)$ is a $(1,1,2,2)$-coloring of $G$.

Conversely, suppose that we have a $(1,1,2,2)$-coloring of $G$ with color classes $W_1, W_2, W_3, W_4$. Since $W_3$ and $W_4$ are $2$-packings in $G$,
$W_3$ and $W_4$ are independent sets in $G^2$. It follows that $W_3 \cup W_4$ induces a bipartite graph in $G^2$.  Let $V_1=W_3 \cup W_4$, $V_2=W_2$, and
$V_3=W_1$.  By definition, $\{V_1, V_2, V_3\}$ is a partition of $V(G)$ as claimed in the statement of the lemma.
\qed

\section{Generalized prisms and the Petersen graph}
\label{sec:generalized-prisms}

In this section we confirm Conjecture~\ref{con:main} for all generalized prisms of cycles, where a {\em generalized prism} is a cubic graph obtained from the disjoint union of two cycles of equal length by adding a perfect matching between the vertices of the two cycles. Along the way we prove that a generalized prism of a cycle is $(1,1,2,2)$-colorable unless it is the Petersen graph, thus characterizing the Petersen graph $P$ in a new way. By separately verifying that $\pch(S(P)) = 5$, Conjecture~\ref{con:main} for generalized prisms then follows from Proposition~\ref{prop:necessary}. We begin with the following technical lemma.

\begin{lemma}\label{lem:help}
If $C_n = v_1\cdots v_n$ is a cycle on $n$ vertices, then the following hold.
\begin{enumerate}
\item[(i)] There exists a set $A \subset V(C_n)$ such that at most one pair of adjacent vertices in $C_n$ is in the complement of $A$ and $G^2[A]$ is an even cycle or a path.
\item[(ii)] If $n$ is odd and $i \in \{3, \dots, n-1\}$, there exists a set $A \subset V(G)$ such that $\{v_1, v_i, v_{j}\} \cap A = \emptyset$ for some $j \in \{i-1, i+1\}$, $v_iv_{j}$ is the only adjacent pair of vertices in $C_n$, which is in the complement of $A$, and $G^2[A]$ is a path.
\end{enumerate}
\end{lemma}

\proof
The result is trivial if $3 \le n \le 5$ so we may assume that $n \ge 6$. To prove statement (i), we first assume $n$ is even, and let
\[A_1 = \begin{cases} \{v_i\mid i \textrm{ is odd}\} & \text{ if } n \equiv 0 \pmod{4}\\ \{v_1, v_2, v_4, v_5\} \cup \{v_j \mid j \text{ is odd and } j \ge 7\} & \text{ if }n\equiv 2\pmod{4}
\end{cases}.\]
Note that $G^2[A_1]$ is an even cycle.

Suppose next that $n$ is odd, $n\ge 7$. If we let
\[A_2 = \begin{cases} \{v_1, v_2, v_4, v_5, v_7, v_8\} \cup \{v_j\mid j \text{ is even and } j \ge 10\} & \text{ if }n\equiv 1 \pmod{4}, n \ge 9\\
\{v_1, v_4\} \cup \{v_j \mid j \text{ is even and }j \ge 4\} & \text{ if } n \equiv 3 \pmod{4}
\end{cases},\]
then $G^2[A_2]$ is a path if $n \equiv 3\pmod{4}$ and $G^2[A_2]$ is an even cycle if $n \equiv 1 \pmod{4}$. This concludes the proof of (i).

We next prove (ii) in which case $n$ is odd.
Let $i \in \{3, \dots, n-1\}$. Suppose first that $i$ is even. If $i \le n-3$, let
\[A_3 = \{v_2, v_{i+2}\} \cup \{v_j\mid 3 \le j \le i-1, j \text{ odd}\} \cup \{v_j \mid i +3 \le j \le n, j \text{ odd}\} , \]  and if $i = n-1$, let
\[ A_4 = \begin{cases} \{v_2, v_{n-3}, v_n\} \cup \{v_j \mid 3 \le j \le n-4, j \text{ odd}\} & \text{ if } n\equiv 1 \pmod{4}\\
\{v_n\} \cup \{v_j\mid 2 \le j \le n-3, j\text{ even}\} & \text{ if } n \equiv 3 \pmod{4}
\end{cases}.\]

Finally, if $i$ is odd, we let $A_5 = \{v_j \mid 2 \le j \le i-1, j \text{ even}\} \cup \{v_j \mid i+2 \le j \le n, j \text{ odd}\}$. In each case, $G^2[A_j]$ is a path for $j \in \{3, 4, 5\}$.
\qed

\begin{theorem}\label{thm:onecycle}
If $G$ is a generalized prism of a cycle, then $G$ is $(1,1,2,2)$-colorable if and only if $G$ is not the Petersen graph.
\end{theorem}

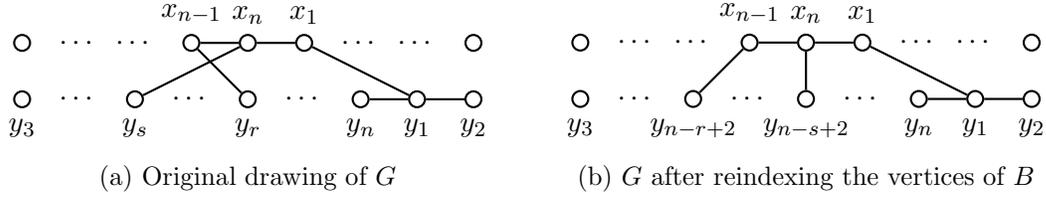
\begin{figure}[ht!]
\begin{center}
\begin{subfigure}[]{.5\textwidth}
\centering
\begin{tikzpicture}[scale=0.75,style=thick]
	\vertex (1) at (0,0) [label=below:$y_3$] {};
	\node(dot1) at (1,0)[]{$\cdots$};
	\vertex (2) at (2,0) [label=below:$y_s$] {};
	\node(dot5) at (3,0)[]{$\cdots$};
	\vertex (3) at (4,0) [label=below:$y_r$] {};
	\node(dot6) at (5,0)[]{$\cdots$};
	\vertex(4) at (6,0)[label=below:$y_n$]{};
	\vertex (5) at (7,0) [label=below:$y_1$] {};
	\vertex (6) at (8,0) [label=below:$y_2$] {};
	\vertex (7) at (0,1) [] {};
	\node(dot3) at (1,1)[]{$\cdots$};
	\node(dot7) at (2,1)[]{$\cdots$};
	\vertex (8) at (3,1) [label=above:$x_{n-1}$] {};
	\vertex (9) at (4,1) [label=above:$x_n$] {};
	\vertex (10) at (5,1) [label=above:$x_1$] {};
	\node(dot4) at (6,1)[]{$\cdots$};
	\node(dot8) at (7,1)[]{$\cdots$};
	\vertex (11) at (8,1) [] {};
	\path
		(8) edge (9)
		(9) edge (10)
		(2) edge (9)
		(3) edge (8)
		(5) edge (10)
		(4) edge (5)
		(5) edge (6)

	;
\end{tikzpicture}
\caption{Original drawing of $G$}
\end{subfigure}~
\begin{subfigure}[]{.5\textwidth}
\centering
\begin{tikzpicture}[scale=0.75,style=thick]
	\vertex (1) at (0,0) [label=below:$y_3$] {};
	\node(dot1) at (1,0)[]{$\cdots$};
	\vertex (2) at (2,0) [label=below:$y_{n-r+2}$] {};
	\node(dot5) at (3,0)[]{$\cdots$};
	\vertex (3) at (4,0) [label=below:$y_{n-s+2}$] {};
	\node(dot6) at (5,0)[]{$\cdots$};
	\vertex(4) at (6,0)[label=below:$y_n$]{};
	\vertex (5) at (7,0) [label=below:$y_1$] {};
	\vertex (6) at (8,0) [label=below:$y_2$] {};
	\vertex (7) at (0,1) [] {};
	\node(dot3) at (1,1)[]{$\cdots$};
	\node(dot7) at (2,1)[]{$\cdots$};
	\vertex (8) at (3,1) [label=above:$x_{n-1}$] {};
	\vertex (9) at (4,1) [label=above:$x_n$] {};
	\vertex (10) at (5,1) [label=above:$x_1$] {};
	\node(dot4) at (6,1)[]{$\cdots$};
	\node(dot8) at (7,1)[]{$\cdots$};
	\vertex (11) at (8,1) [] {};
	\path
		(8) edge (9)
		(9) edge (10)
		(2) edge (8)
		(3) edge (9)
		(5) edge (10)
		(4) edge (5)
		(5) edge (6)

	;
\end{tikzpicture}
\caption{$G$ after reindexing the vertices of $B$}
\end{subfigure}
\caption{Drawing of the generalized prism of a cycle}
\label{fig:drawing}
\end{center}
\end{figure}

\proof
Up to isomorphism there is only one generalized prism of the $3$-cycle, and it is clearly $(1,1,2,2)$-colorable. So we may assume that $C_n$ is a cycle
on at least four vertices. By Lemma~\ref{lem:CP1}, it suffices to show that $V(G)$ can be partitioned into $V_1, V_2, V_3$, where $V_2$ and $V_3$ are independent sets and $G^2[V_1]$ is bipartite. In $G$, let $A = x_1 \cdots x_n$ and $B = y_1 \cdots y_n$ represent the two copies of $C_n$. By definition, there exists a perfect matching between $A$ and $B$ in $G$, and so we define $f:A \to B$ such that $f(x_i) = y_j$ if $x_iy_j \in E(G)$ for all $1 \le i \le n$. Without loss of generality we assume that $f(x_1) = y_1$.  In addition,
let $f(x_{n-1}) = y_r$ and $f(x_n) = y_s$ for some $\{r, s\} \subseteq \{2, \dots, n\}$. We then draw $A$ horizontally so that the indices increase from left to right and $x_n$ is located in the middle of the cycle. Moreover, we can draw $B$ horizontally and beneath $A$ so that the indices increase from left to right
and $y_1$ is drawn to the right of both $y_r$ and $y_s$, as shown in Figure~\ref{fig:drawing}(a). If $y_r$ is to the right of $y_s$, meaning $r >s$, then we can relabel the vertices of $B$ so that $f(x_1)$ still has index $1$, but the indices increase from right to left. Then we draw $B$ so that the indices increase from left to right, as depicted in Figure~\ref{fig:drawing}(b), and $y_r$ is to the left of both $y_1$ and $y_s$. So we may assume throughout the remainder of the proof that $r < s$.

Suppose first that $C_n$ is an even cycle. Let
\begin{itemize}
\item $X_2 = \{x_i \mid 1 \le i \le n,\, i \text{ is odd}\}$ and
\item $X_3 = \{x_i \mid 1 \le i \le n, \, i \text{ is even}\}$.
\end{itemize}
Let $Y_1$ be the set $A_1$ from Lemma~\ref{lem:help} and define
\begin{itemize}
\item $Y_2 = \{y_i \not\in Y_1 \mid  f^{-1}(y_i) \in X_3\}$ and
\item $Y_3 =   \{y_i \not\in Y_1 \mid f^{-1}(y_i) \in X_2\}$.
\end{itemize}
One can easily verify that $G^2[Y_1]$ is an even cycle and $X_i \cup Y_i$ is independent for $i \in \{2, 3\}$. By Lemma~\ref{lem:CP1} it
follows that the generalized prism is $(1,1,2,2)$-colorable.

So we may assume that $C_n$ is an odd cycle. Let $X_1 = \{x_n\}$,
\begin{itemize}
\item $X_2 = \{x_i \mid i\text{ is odd}, 1 \le i\le n-2\}$, and
\item $X_3 = \{x_i \mid i\text{ is even}, 1 \le i \le n-1\}$.
\end{itemize}
 In what follows, we partition the vertices of $B$ into $Y_1, Y_2, Y_3$ depending on the parity of $r$ and $s$. In each case, we let $V_i = X_i \cup Y_i$ for each $i \in [3]$ so that $V_2$ and $V_3$ are independent and $G^2[V_1]$ is bipartite.

\noindent{\bf Case $1$.} Suppose $s$ is odd.\\
We shall assume first that $r$ is odd as well. If $s \ne n$, let  $Y_1$ be the set $A_5$ from Lemma~\ref{lem:help} where $i = s$ so that $G^2[V_1]$ is a path that does not contain $\{y_1, y_s, y_{s+1}\}$. We then define

\[Y_2 = \begin{cases} \{y_{s+1}\} \cup \{y_i \notin Y_1 \mid  f^{-1}(y_i)\in X_3\} & \text{ if }f^{-1}(y_{s+1}) \in X_3\\
\{y_s\} \cup \{y_i \not\in Y_1 \mid  f^{-1}(y_i) \in X_3\} & \text{ if } f^{-1}(y_{s+1}) \in X_2
\end{cases}, \]
and

\[Y_3 = \begin{cases} \{y_s\} \cup \{y_i \not\in Y_1 \mid   f^{-1}(y_i) \in X_2\} & \text{ if }f^{-1}(y_{s+1}) \in X_3\\
\{y_{s+1}\} \cup \{y_i \not\in Y_1 \mid   f^{-1}(y_i) \in X_2\} & \text{ if } f^{-1}(y_{s+1}) \in X_2
\end{cases}.\]
If $s = n$ and $r$ is odd or $r = n-1$, we let $Y_1 = \{y_i \mid  i \text{ is even}\}$,

\begin{itemize}
\item $Y_2 = \{y_s\} \cup \{y_i \not\in Y_1 \mid  f^{-1}(y_i) \in X_3\}$ and
\item $Y_3 = \{y_i \not\in Y_1 \mid  f^{-1}(y_i) \in X_2\}$.
\end{itemize} 
In either case, $G^2[V_1]$ is a path. So we may assume that $r$ is even and $r \ne n-1$ (note that the case $r=n-1$ is symmetric to the case when $s=n$, which was considered above).

 If $s < n$, we let

 \begin{itemize}
\item $Y_1 = \{y_i \mid 3 \le i \le s-2, i \text{ is odd}\} \cup \{y_2, y_{s+1}\} \cup \{y_i \mid s+2 \le i \le n, i \text{ is odd}\}$.
 \end{itemize}
 As above, $G^2[V_1]$ is a path. If $f^{-1}(y_{s-1}) \in X_3$, then let
 \begin{itemize}
 \item $Y_2 = \{y_i \not\in Y_1 \mid  f^{-1}(y_i) \in X_3\}$ and
 \item $Y_3 = \{y_s\} \cup \{y_i \not\in Y_1 \mid  f^{-1}(y_i) \in X_2\}$.
 \end{itemize}
Otherwise, let
\begin{itemize}
\item $Y_2 = \{y_s\} \cup \{y_i \not\in Y_1 \mid  f^{-1}(y_i) \in X_3\}$ and
\item $Y_3 = \{y_i \not\in Y_1 \mid  f^{-1}(y_i) \in X_2\}$.
\end{itemize}

\begin{figure}[ht!]
\begin{center}
\begin{subfigure}[]{.5\textwidth}
\centering
\begin{tikzpicture}[scale=0.75,style=thick]
	\vertex (1) at (0,0) [] {};
	\vertex (2) at (1,0) [label=below:1] {};
	\vertex (3) at (2,0) [label=below:1] {};
	\vertex(4) at (3,0)[]{};
	\vertex (5) at (4,0) [] {};
	\vertex (6) at (5,0) [label=below:1] {};
	\vertex (7) at (6,0) [] {};
	\vertex (8) at (7,0) [label=below:1] {};
	\vertex (9) at (8,0) [label=below:1] {};
	\node(dot1) at (2,1)[]{$\cdots$};
	\vertex(10) at (3,1)[label=above:$x_{n-1}$]{};
	\vertex(11) at (4,1)[label=above:$x_n$]{};
	\vertex(12) at (5,1)[label=above:$x_1$]{};
	\node(dot2) at (6,1)[]{$\cdots$};
	\path
		(1) edge (2)
		(2) edge (3)
		(3) edge (4)
		(4) edge (5)
		(5) edge (6)
		(6) edge (7)
		(7) edge (8)
		(8) edge (9)
		(1) edge[bend right=20] (9)
		(10) edge (11)
		(11) edge (12)
		(10) edge (3)
		(11) edge (5)
		(12) edge (7)

	;
\end{tikzpicture}
\caption{$Y_1$ when $n=9$}
\end{subfigure}~
\begin{subfigure}[]{.5\textwidth}
\centering
\begin{tikzpicture}[scale=0.75,style=thick]
	\vertex (1) at (0,0) [] {};
	\vertex (2) at (1,0) [label=below:1] {};
	\vertex (3) at (2,0) [] {};
	\vertex(4) at (3,0)[label=below:1]{};
	\vertex (5) at (4,0) [] {};
	\vertex (6) at (5,0) [] {};
	\vertex (7) at (6,0) [label=below:1] {};
	\vertex (8) at (7,0) [] {};
	\vertex (9) at (8,0) [label=below:1] {};
	\vertex(10) at (9,0)[]{};
	\vertex(11) at (10,0)[label=below:1]{};
	\node(dot1) at (3,1)[]{$\cdots$};
	\vertex(12) at (4,1)[label=above:$x_{n-1}$]{};
	\vertex(13) at (5,1)[label=above:$x_n$]{};
	\vertex(14) at (6,1)[label=above:$x_1$]{};
	\node(dot2) at (7,1)[]{$\cdots$};
	\path
		(1) edge (2)
		(2) edge (3)
		(3) edge (4)
		(4) edge (5)
		(5) edge (6)
		(6) edge (7)
		(7) edge (8)
		(8) edge (9)
		(9) edge (10)
		(10) edge (11)
		(1) edge[bend right=20] (11)
		(12) edge (13)
		(13) edge (14)
		(12) edge (4)
		(13) edge (6)
		(14) edge (8)

	;
\end{tikzpicture}
\caption{$Y_1$ when $n=11$}
\end{subfigure}
\caption{The set $Y_1$ in the (1,1,2,2)-coloring of $C_9$ and $C_{11}$}
\label{fig:boat}
\end{center}
\end{figure}
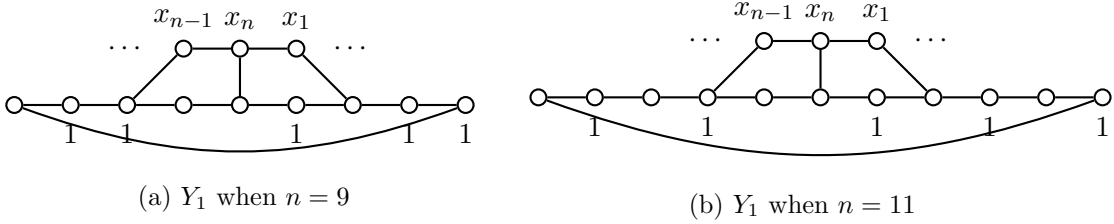

\noindent{\bf Case $2$.} Suppose that $s$ is even.\\
First, note that if $r$ is odd, then we can define $Y_1, Y_2, Y_3$ similarly to those subcases given in Case 1 (which can be observed by reversing the roles of $x_{n-1}$ and $x_1$). So we may assume that $r$ is even. Suppose first that $n >5$ and $r>2$. Then one of the sets $A_3$ or $A_4$ given in Lemma~\ref{lem:help} can be chosen for $Y_1$ so that $G^2[V_1]$ is a path and $Y_1$ does not contain vertices $\{y_1, y_s, y_j\}$ where $j \in \{s-1, s+1\}$ (see Figure~\ref{fig:boat} for two corresponding examples). Thus, $G^2[V_1]$ is path. We then define
\[Y_2 = \begin{cases} \{y_i\not\in Y_1 \mid  f^{-1}(y_i) \in X_3\} & \text{ if }s \le n-3, f^{-1}(y_{s+1}) \in X_3\\

\{y_s\} \cup \{y_i \not\in Y_1\mid  f^{-1}(y_i) \in X_3\} & \text{ if } s \le n-3, f^{-1}(y_{s+1}) \in X_2\\

\{y_i \not\in Y_1 \mid  f^{-1}(y_i) \in X_3\} & \text{ if } s = n-1, f^{-1}(y_{s-1}) \in X_3\\

\{y_s\} \cup  \{y_i \not\in Y_1 \mid  f^{-1}(y_i) \in X_3\} & \text{ if } s = n-1, f^{-1}(y_{s-1}) \in X_2
\end{cases},\]
 and

 \[Y_3  = \begin{cases} \{y_s\} \cup \{y_i \not\in Y_1 \mid  f^{-1}(y_i) \in X_2\} & \text{ if } s \le n-3, f^{-1}(y_{s+1}) \in X_3\\

 \{y_i \not\in Y_1 \mid  f^{-1}(y_i) \in X_2\} & \text{ if } s \le n-3, f^{-1}(y_{s+1}) \in X_2\\

 \{y_s\} \cup \{y_i \not\in Y_1 \mid  f^{-1}(y_i) \in X_2\} & \text{ if } s = n-1, f^{-1}(y_{s-1}) \in X_3\\

  \{y_i \not\in Y_1 \mid  f^{-1}(y_i) \in X_2\} & \text{ if } s=n-1, f^{-1}(y_{s-1}) \in X_2
 \end{cases}.\]

Next, we suppose that $r=2$ (while we may assume, by symmetry, that $s<n-1$), and let
\[Y_1 = \{y_{s+2}\} \cup \{y_i \mid s+3 \le i \le n, i \text{ odd}\} \cup \{y_i \mid 3 \le i \le s-1, i \text{ odd}\}.\]
Then choose $Y_2$ and $Y_3$ in the same way as above based on the index of $s$.

\begin{figure}[h!]
\begin{center}
\begin{tikzpicture}[]
\tikzstyle{vertex}=[circle, draw, inner sep=0pt, minimum size=6pt]
\tikzset{vertexStyle/.append style={rectangle}}
	\vertex (1) at (0,0) [label=below:$2$] {};
	\vertex (2) at (1,0) [label=below:$3$] {};
	\vertex (3) at (2,0) [label=below:$1$] {};
	\vertex (4) at (3,0) [label=below:$2$] {};
	\vertex (5) at (4,0) [label=below:$3$] {};
	\vertex (6) at (0,1) [label=above:$2$] {};
	\vertex (7) at (1,1) [label=above:$3$] {};
	\vertex (8) at (2,1) [label=above:$1$] {};
	\vertex (9) at (3,1) [label=above:$2$] {};
	\vertex (10) at (4,1) [label=above:$3$] {};
	\path
		(1) edge (2)
		(2) edge (3)
		(3) edge (4)
		(4) edge (5)
		(1) edge[bend right=20] (5)
		(6) edge (7)
		(7) edge (8)
		(8) edge (9)
		(9) edge (10)
		(6) edge[bend left=20] (10)
		(1) edge (7)
		(2) edge (6)
		(3) edge (8)
		(4) edge (10)
		(5) edge (9)

	;
\end{tikzpicture}
\end{center}
\caption{The labels depict a $(1,1,2,2)$-coloring}
\label{fig:5cyclelabel}
\end{figure}
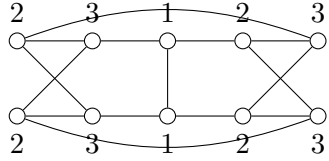

Finally, consider when $C_n$ is a $5$-cycle. If $f(x_2) = y_5$ and $f(x_3) = y_3$, then Figure~\ref{fig:5cyclelabel} depicts a
labeling of $G$ where $V_2$ and $V_3$ are independent and $G^2[V_1]$ is bipartite. If $f(x_2) = y_3$ and $f(x_3) = y_5$, then $G$ is the Petersen graph. The argument is complete by invoking the fact~\cite[Proposition 4]{gt-2016} that the Petersen graph is not $(1,1,2,2)$-colorable.
\qed

\begin{figure}[h!]
\begin{center}
\begin{tikzpicture}[scale=1.05]
\tikzstyle{vertex}=[circle, draw, inner sep=0pt, minimum size=.45cm]
\tikzset{vertexStyle/.append style={rectangle}}
	\vertex(1) at (-.15,0)[]{2};
	\vertex(2) at (2,0)[]{1};
	\vertex(3) at (4.15,0)[]{4};
	\vertex(4) at (-1.5, 4)[]{3};
	\vertex(5) at (5.5, 4)[]{1};
	\vertex(6) at (2, 6.8)[]{1};
	\vertex(7) at (2, 4.95)[]{1};
	\vertex(8) at (.3, 3.4)[]{2};
	\vertex(9) at (3.7, 3.4)[]{1};
	\vertex(10) at (.9, 1.2)[]{1};
	\vertex(11) at (3.1, 1.2)[]{1};
	\vertex(12) at (0.15, 5.33)[]{2};
	\vertex(13) at (3.85,5.33)[]{5};
	\vertex(14) at (-.825,2)[]{1};
	\vertex(15) at (4.83,2)[]{2};
	\vertex(16) at (.375, .6)[]{5};
	\vertex(17) at (3.625,.6)[]{3};
	\vertex(18) at (2, 5.8)[]{4};
	\vertex(19) at (-.52,3.68)[]{1};
	\vertex(20) at (4.52, 3.68)[]{3};
	\vertex(21) at (2.35,2.45)[]{2};
	\vertex(22) at (1.65,2.45)[]{5};
	\vertex(23) at (1.45,3)[]{3};
	\vertex(24) at (2.55,3)[]{2};
	\vertex(25) at (2,3.4)[]{4};
	
	\path
		(1) edge (2)
		(2) edge (3)
		(1) edge (14)
		(14) edge (4)
		(3) edge (15)
		(15) edge (5)
		(4) edge (12)
		(12) edge (6)
		(5) edge (13)
		(13) edge (6)
		(6) edge (18)
		(18) edge (7)
		(4) edge (19)
		(19) edge (8)
		(5) edge (20)
		(20) edge (9)
		(1) edge (16)
		(16) edge (10)
		(3) edge (17)
		(17) edge (11)
		(7) edge (23)
		(23) edge (10)
		(10) edge (21)
		(21) edge (9)
		(11) edge (22)
		(22) edge (8)
		(8) edge (25)
		(9) edge (25)
		(7) edge (24)
		(11) edge (24)

	;
\end{tikzpicture}
\end{center}
\caption{$5$-packing coloring of $S(P)$}
\label{fig:Petersen}
\end{figure}
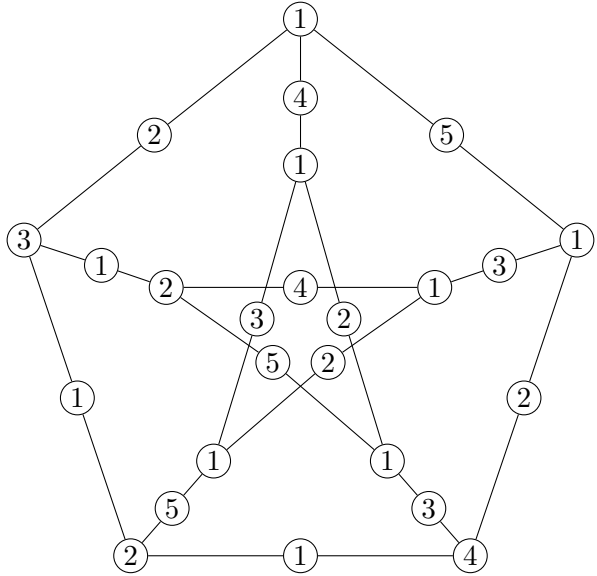

By Theorem~\ref{thm:onecycle} and Proposition~\ref{prop:necessary} we know that any subdivided generalized prism of a cycle but the subdivided Petersen graph $P$ is $5$-packing colorable. In addition, a $5$-packing coloring of $S(P)$ is shown in Fig.~\ref{fig:Petersen}. Hence we have the following result.

\begin{corollary}
If $G$ is a generalized prism of a cycle, then $\pch(S(G))\le 5$.
\end{corollary}

Intuitively, it seems reasonable to expect that an optimal packing coloring of any subdivided graph colors all the subdivided vertices by $1$. The example from Fig.~\ref{fig:Petersen} shows that an optimal coloring need not be like that. In fact, no optimal coloring of $S(P)$ colors all the subdivided vertices by $1$. This is an immediate consequence of the following result.

\begin{proposition}\label{prp:uglycoloring}
If $G$ is not $(1,1,2,2)$-colorable, and $S(G)$ is $(1,2,3,4,5)$-colorable, then in every $5$-packing coloring of $S(G)$ at least
one of the subdivided vertices of $S(G)$ receives color bigger than $1$.
\end{proposition}
\proof
Suppose to the contrary that $c$ is a $5$-packing coloring of $S(G)$ with $c(u_e)=1$ for every edge $e\in E(G)$. Then all vertices of $V(G)$ in $S(G)$ receive colors from $\{2,3,4,5\}$. Consider the coloring $c'$ of $V(G)$, obtained as the restriction of $c$ to $G$. Note that vertices colored by the color $2$, respectively $3$, form an independent set in $G$, while the set of vertices colored by the color $4$, respectively $5$, is a 2-packing of $G$. This implies that $c'$ is a $(1,1,2,2)$-coloring of $G$.
\qed

\section{A class larger than generalized prisms}
\label{sec:more-general}

In this section we confirm Conjecture~\ref{con:main} for a class of graphs larger than generalized prisms of cycles by proving that if $G$ is a connected, cubic graph of order $2n$ with a $2$-factor $\mathcal{F}$ and a perfect matching $M$, where $\mathcal{F}$ contains a cycle $C$ of length $n$, no edge of $M$ has both vertices in $C$, and $\mathcal{F}$ contains at most one $5$-cycle, then $G$ is $(1,1,2,2)$-colorable. In other words, our result extends Theorem~\ref{thm:onecycle} to the graphs obtained from generalized prisms in such a way that one of the two $n$-cycles in the edge set of a generalized prism is replaced by a union of cycles among which at most one is a 5-cycle.

\begin{theorem} \label{thm:2connectedcubic}
Let $G$ be a connected, cubic graph of order $2n$ with a $2$-factor $\mathcal{F}$ and a perfect matching $M$.
If $\mathcal{F}$ contains a cycle $C$ of length $n$ where no edge of $M$ has both vertices in $C$, and $\mathcal{F}$ contains at
most one $5$-cycle, then $G$ is $(1,1,2,2)$-colorable.
\end{theorem}

\proof
Note that by Theorem~\ref{thm:onecycle}, we may assume that $\mathcal{F}$ contains at least three cycles.
Thus, $n \ge 6$. We let $C_n=x_1\cdots x_n$ represent the cycle in $\mathcal{F}$ of order $n$, and let $Z_1, \dots, Z_k$ be
the remaining cycles of $\mathcal{F}$.  Reindexing if necessary, we may assume that if $\mathcal{F}$ contains a $5$-cycle that $Z_1$ is said $5$-cycle. Otherwise, if $Z_i$ is odd for some $i \in [k]$, we let $Z_1$ represent the smallest odd cycle among all $Z_i$, $i \in [k]$. In any case, we let $Z_1 = y_1\cdots y_p$ for some $3 \le p <n$.

Assume first that $Z_1$ is a $5$-cycle so that $p = 5$. Note that there exists $x_i \in C_n$ such that
$f(x_i) \in Z_1$ and $f(x_j) \not\in Z_1$ for some $j \in \{i-1, i+1\}$. Reindexing $x_1, \dots, x_n$ if necessary, we may
assume $f(x_n) \in Z_1$ and, redrawing $G$ if necessary, $f(x_{n-1}) \not\in Z_1$.

As in Theorem~\ref{thm:onecycle}, we let $X_1 = \{x_n\}$, $X_2 = \{x_i \mid 1 \le i < n, i \text{ is odd}\}$ and $X_3 = \{x_i \mid 2 \le i <n, i \text{ is even}\}$. In what follows, we partition the vertices of $\cup_{i=1}^k Z_i$ into $Y_1, Y_2, Y_3$, and let $V_i = X_i \cup Y_i$ for each $i \in [3]$. In each case, $G^2[V_1]$ will be bipartite and $V_2, V_3$ will be independent sets.

\noindent{\bf Case 1.} Suppose that $f(x_1) \in Z_1$. Without loss of generality, we may assume $f(x_1) = y_1$, and
reindexing $Z_1$ if necessary, $f(x_n) = y_s$ where $s \in \{4,5\}$. For each $i \in \{2, \dots, k\}$ let $T_i$ be one of the sets $A_1$ or $A_2$ from Lemma~\ref{lem:help} depending on the congruence class of $n$ modulo $4$. Note that for each $i \in \{2, \dots, k\}$, $G^2[T_i]$ is bipartite.

Next, we assume for the time being that $s = 4$ and let $T_1 = \{ y_2, y_5\}$. We let
\begin{itemize}
\item $Y_1 = \cup_{i=1}^k T_i$,
\item $Y_2 = W_2 \cup \bigcup_{i=2}^k \{v_j \in Z_i - T_i \mid f^{-1}(v_j) \in X_3\}$ and
\item $Y_3 = W_3 \cup \bigcup_{i=2}^k \{v_j \in Z_i - T_i \mid f^{-1}(v_j) \in X_2\}$,
\end{itemize}
where
\begin{itemize}
\item $W_2 = \{y_3\}$ and $W_3 = \{y_1,y_4\}$, if $f^{-1}(y_3) \in X_3$;
\item $W_2 = \{y_4\}$ and $W_3 = \{y_1, y_3\}$, if $f^{-1}(y_3) \in X_2$.
\end{itemize}

In $G^2[V_1]$, all edges incident to $x_n$ are bridges to either the $K_2$ induced by $T_1 - \{x_n\}$ or to a
bipartite component induced by $T_i$ for some $i \in \{2, \dots, k\}$. Thus, $G^2[V_1]$ is bipartite.
Furthermore, $V_i$ where $i \in \{2, 3\}$ is independent.

Now, one can easily see that $Y_1, Y_2, Y_3$ can be defined in a similar fashion if instead $s = 5$. \smallqed

\noindent{\bf Case 2.} Suppose that $f(x_1) \not\in Z_1$. Without loss of generality, we may assume $f(x_n) = y_1$. In this case,
we define $T_i$ as in Case 1 for each $i \in \{2, \dots, k\}$ and we let $T_1 = \{y_2, y_4\}$. We let
\begin{itemize}
\item $Y_1 = \bigcup_{i=1}^k T_i$,
\item $Y_2 = W_2 \cup \bigcup_{i=2}^k \{v_j \in Z_i - T_i \mid f^{-1}(v_j) \in X_3\}$ and
\item $Y_3 = W_3 \cup \bigcup_{i=2}^k \{v_j \in Z_i - T_i \mid f^{-1}(v_j) \in X_2\}$,
\end{itemize}
where $y_3\in W_2$ if and only if $f^{-1}(y_3) \in X_3$, and otherwise $y_3\in W_3$; and $y_1$ and $y_5$ are in different sets $W_2,W_3$, depending on $f^{-1}(y_5)$.

As in Case 1, $G^2[V_1]$ is bipartite and $V_i$ is independent for $i \in \{2,3\}$. \smallqed
\smallskip

Now consider the case, when at least one of the cycles $Z_i$ is odd and none of them is a $5$-cycle. Recall that $Z_1$ is a shortest odd cycle from $\mathcal{F}$. We shall assume that $f(x_n) = y_s$ for some $s \in [p]$. Whether or not $f(x_1) \in Z_1$, we may choose $T_1$ to be the set $A_3, A_4$, or $A_5$ from Lemma~\ref{lem:help} so that $\{y_s, y_{s+1}\}\cap T_1 = \emptyset$ if $f(x_1) \not\in Z_1$, $\{f(x_1), y_s, y_{s+1}\} \cap T_1 = \emptyset$ if $f(x_1) \in Z_1$, and $G^2[T_1 \cup X_1]$ is bipartite. Then for each $i \in \{2, \dots, k\}$, we let $T_i$ be one of the sets $A_1$ or $A_2$ from Lemma~\ref{lem:help} depending on the congruence class of $n$ modulo $4$. Defining $Y_1$, $Y_2$, and $Y_3$ similarly as in Case 1, one can verify that $G^2[V_1]$ is indeed bipartite.

Finally, consider the case that $Z_1$ is even, in which case all the cycles $Z_i$ are even, and so $n$ is also even. In this and only in this case, we let $X_1=\emptyset$, and $X_2 = \{x_i \mid 1 \le i \le n-1, i \text{ is odd}\}$ and $X_3 = \{x_i \mid 2 \le i \le n, i \text{ is even}\}$. Next, for each $i \in [k]$, we let $T_i$ be the set $A_1$ from Lemma~\ref{lem:help} and we define $Y_1 =  \bigcup_{i=1}^k T_i$ (note that $V_1=Y_1$).
Letting
\begin{itemize}
\item $Y_2 = \bigcup_{i=1}^k \{v_j \in Z_i - T_i \mid f^{-1}(v_j) \in X_3\}$ and
\item $Y_3 = \bigcup_{i=1}^k \{v_j \in Z_i - T_i \mid f^{-1}(v_j) \in X_2\}$
\end{itemize}
we obtain a $(1,1,2,2)$-coloring of $G$.
\qed

Note that the graphs from Theorem~\ref{thm:2connectedcubic} are $2$-connected. We suspect that a similar approach might work to prove that an arbitrary 2-connected cubic graph (except the Petersen graph) has a $(1,1,2,2)$-packing coloring.
(Recall that by Petersen's theorem~\cite{petersen} the edge set of any such graph can be partitioned into a 2-factor and a perfect matching.)

One class of cubic graphs covered by the result in Theorem~\ref{thm:2connectedcubic} are some subclasses of generalized Petersen graphs.
Let $k$ and $n$ be positive integers such that $k < n/2$.  The \emph{generalized Petersen graph} $P(n,k)$ has vertex set
$\{u_1,v_1, \ldots, u_n,v_n\}$.  The edge set of $P(n,k)$ is the set
\[ \{u_iu_{i+1} \mid i \in [n]\} \cup \{u_iv_i \mid i \in [n]\} \cup \{v_iv_{i+k} \mid i \in [n]\}\,,\]
where addition on the subscripts is computed modulo $n$.  The set $\{u_i \mid i\in [n]\}$ induces a cycle of order $n$,
while the set $\{v_i \mid i \in [n]\}$ induces a disjoint union of cycles.  The order and the number of this latter collection
of cycles depends on the relationship between $n$ and $k$.  It is easy to see that if $n$ and $k$ are relatively prime,
then $\{v_i \mid i \in [n]\}$ induces a single cycle of order $n$.  In this case $P(n,k)$ is a generalized prism of $C_n$ and satisfies the hypotheses of
Theorem~\ref{thm:onecycle} unless $n=5$, in which case $P(5,k)$ is either the ordinary prism of $C_5$ (that is, the
Cartesian product of $C_5$ and $K_2$) or the famous Petersen graph.  If $n$ and $k$ are not relatively prime, then the subgraph
of $P(n,k)$ induced by $\{v_i \mid i \in [n]\}$ consists of the disjoint union of $n/r$ cycles each of order $r$, where
$r$ is the smallest positive integer such that $rk$ is divisible by $n$.  Hence, these will be $5$-cycles if and only if
$n$ is a multiple of $5$.
\begin{corollary} \label{cor:GPG}
If $n$ and $k$ are positive integers such that $k<n/2$ and $n$ is not a multiple of $5$, then $P(n,k)$ has a $(1,1,2,2)$-coloring and
hence $\pch(S(P(n,k)))\le 5$.
\end{corollary}

\section{Multiple subdivisions}
\label{sec:multiple}

We have already remarked that $\pch(S(K_n))=n+1$. We next consider $\pch(S_i(K_n))$ for $i\ge 2$.

\begin{proposition} \label{prop:SubofKn}
If $n\ge 3$ and $i\ge 3$, then
$$\pch(S_i(K_n))=\left\{
  \begin{array}{ll}
    3; & \hbox{if $i\equiv 3 \pmod 4  $,} \\
    4; & \hbox{otherwise.}
  \end{array}
\right.$$
Moreover, $\pch(S_2(K_n))\xrightarrow[n\to \infty]{}\infty$.
\end{proposition}

\proof
Clearly, $\pch(S_i(K_n))\ge 3$ for $n\ge 3$ and $i\ge 3$. Note that $S_i(K_n)$ contains a cycle of length $3i + 3 = 3(i+1)$.
Since $\pch(C_n) = 3$ if $n\equiv 0 \pmod 4$, and $\pch(C_n) = 4$ otherwise (see~\cite{goddard-2008}), we get that $\pch(S_i(K_n))\ge 3$ if $i\equiv 3\pmod 4$, and $\pch(S_i(K_n))\ge 4$ otherwise.

To prove that these lower bounds are tight, we color $S_i(K_n)$ as follows. If $i\equiv 3\pmod 4$, then color the vertices $v\in V(K_n)$ with color $3$;  otherwise color all these vertices with $4$. Colorings of the subdivided vertices are done based on the parity of $i\pmod 4$ as follows. If $i\equiv 3 \pmod 4$, then for each original edge of $K_n$ color the subdivided vertices consecutively by $1,2,1$, and add the block of colors $3,1,2,1$ as many times as required. If $i=4$, use colors $1,2,3,1$. For any even $i\ge 6$, alternatively attach to the four colors $1,2,3,1$ the pairs $2,1$ and $3,1$ as many times as required. Finally, let $i\equiv 1 \pmod 4$. If $i = 5$, then use the pattern $1,3,1,2,1$, and if $i\ge 9$, then add the block $3,1,2,1$ as many times as required. In all of the cases it is straightforward to verify that the constructed colorings are packing colorings.

It remains to consider the case $i=2$. If $e=uv\in E(K_n)$, then let $u_e$ and $v_e$ be the vertices of $S_2(K_n)$ obtained by subdividing the edge $e$, where $u_e$ is adjacent to $u$ and $v_e$ to $v$. Let $c$ be an arbitrary packing coloring of $S_2(K_n)$.  Then for any edge $e=uv\in E(K_n)$ we must have $c(u_e)\ne 1$ or $c(v_e)\ne 1$. Define now the orientation of $K_n$ as follows. If for the edge $e=uv$ we have $c(u_e)\ne 1$, then in $K_n$ orient the edge $uv$ from $u$ to $v$. Otherwise we must have $c(v_e)\ne 1$
 in which case we orient the edge $uv$ from $v$ to $u$. (In the case that both $c(u_e)\ne 1$ and $c(v_e)\ne 1$ hold, we orient the edge $uv$ arbitrarily.) By the degree sum formula for digraphs and the pigeon-hole principle there exists a vertex $u$ with out-degree at least $\lceil (n-1)/2 \rceil$. This means that in $S_2(K_n)$ $u$ has at least that many neighbors colored with different colors bigger than $1$. Hence $\pch(S_2(K_n)) > \lceil (n-1)/2 \rceil$.
\qed

By using our earlier observation that the packing chromatic number of a subgraph is bounded above by the packing
chromatic number of the original graph, we get the following immediate corollary of Proposition~\ref{prop:SubofKn}.

\begin{corollary}
\label{cor:subdivisiongraph}
If $G$ is a connected graph of order at least $3$ and $i \ge 3$, then
\[3 \le \pch(S_i(G)) \le 4\,.\]
\end{corollary}

In the case of trees we can further strengthen the result of  Corollary~\ref{cor:subdivisiongraph} by including the parameter $i=2$ (and $i=1$), and by showing that for any odd $i$, the packing chromatic number is always $3$. More precisely, we have the following result.

\begin{theorem}
\label{thm:subdivided-trees}
If $i\ge 1$, then
$$\max\{\pch(S_i(T))\ |\ T\ {\rm tree}\} = \left\{
  \begin{array}{ll}
    3; & i\ {\rm odd}, \\
    4; & i\ {\rm even}.
  \end{array}
\right.$$
\end{theorem}

\proof
Let $T$ be a tree on at least three vertices. Then, as already mentioned, $\pch(S_1(T))=3$, hence the assertion holds for $i=1$.

Let $i=2$ and let $T$ be an arbitrary tree. To see that $\pch(S_2(T))\le 4$ let $v$ be an arbitrary vertex of $T$ and consider the BFS-tree of $S_2(T)$ rooted in $v$. Then set
$$c(x) = \left\{
  \begin{array}{ll}
    1; & d_T(v,x) \equiv 1 \pmod 3, \\
    2; & d_T(v,x) \equiv 2 \pmod 3, \\
    3; & d_T(v,x) \equiv 0 \pmod 6, \\
    4; & d_T(v,x) \equiv 3 \pmod 6. \\
  \end{array}
\right.$$
It is straightforward to verify that $c$ is a packing coloring of $S_2(T)$. Let now $T$ be a tree with a vertex $u$ of degree at least $3$, let $v$ be a
neighbor of $u$ and let $w$ be a neighbor of $v$ different from $u$.
Recall that if $xy\in E(T)$, then we denote with $e_{xy}$ and $e_{yx}$ the vertices of $S_2(T)$ obtained by subdividing $xy$,
where $e_{xy}$ is the vertex adjacent to $x$. Let $c$ be a packing coloring of $S_2(T)$. If $c(u) = 1$, then considering the neighbors of $u$ (in $S_2(T)$) we see that $\pch(S_2(T))\ge 4$. The same conclusion also follows if $\{c(u), c(v)\} = \{2,3\}$. Suppose next that $c(u) = c(v) = 2$. Then we may without loss of generality assume that $c(e_{vu}) = 3$ which in turn implies that $c(e_{vw}) = 1$. But then $c(e_{wv}) \ge 4$. Finally, let $c(u) = c(v) = 3$. Assuming without loss of generality that $c(e_{vu}) = 1$, we get $c(e_{vw}) = 2$, $c(e_{wv}) = 1$, but then $c(w) \ge 4$. This settles the case $i=2$.

Suppose now that $i\ge 3$. Then by Corollary~\ref{cor:subdivisiongraph}, $\pch(S_i(T))\le 4$. We first deal with $i$ odd in which case we need to prove that $\pch(S_i(T))\le 3$ for any tree $T$. In the first subcase assume that $i\equiv 3\pmod 4$.  By Corollary~\ref{cor:subdivisiongraph} we know
that $3 \le \pch(S_i(T))$.  Since $\pch(S_i(T))\le \pch(S_i(K_n))$, where $T$ has order $n$, and $\pch(S_i(K_n))=3$ by Proposition~\ref{prop:SubofKn},
we conclude that $\pch(S_i(T))=3$ when $i\equiv 3\pmod 4$.
The second subcase to consider is when $i\equiv 1\pmod 4$, $i\ge 5$. Again root $S_i(T)$ in a vertex of $T$, say $u$, and consider the corresponding BFS tree. Consider the following sequence $S$ of $i+1$ colors: first repeat the block $2,1,3,1$ as many times as necessary and finish it with colors $2,1,3$. Note that $|S| \equiv 3 \pmod 4$. Let now $e=xy$ be an edge of $T$, where $d_T(x,u) < d_T(y,u)$. If $d_T(x,u)$ is even, then color the vertices in $S_i(T)$ between $x$ and $y$ (including $x$ and $y$) with the sequence of colors $S$, otherwise (if $d_T(x,u)$ is odd), color the vertices in $S_i(T)$ between $x$ and $y$ (including $x$ and $y$) with the sequence of colors obtained by reversing $S$. Note that this gives a well-defined coloring of $V(S_i(T))$, that is, each vertex of $T$ receives a unique color and that $c$ is a packing coloring.

It remains to consider the case when $i\ge 4$ is even. To complete the argument we need to show that $\pch(S_i(T))\ge 4$ for some tree $T$. Suppose on the contrary that $\pch(S_i(T)) =  3$ holds for any tree $T$ on at least three vertices. If $u$ is a vertex of degree at least $3$, then as above we infer that if $c$ is a $3$-packing coloring of $S_i(T)$, then $c(u) > 1$. In the rest we will also use the fact that if $c(x_1) = 3$ for some vertex of $S_i(T)$, and $x_1, x_2, x_3, x_4, \ldots$  is a path in $S_i(T)$, then $c(x_2) = 1$. Indeed, for otherwise $c(x_2) = 2$, but then $c(x_3) = 1$ and we would have $c(x_4) \ge 4$. Consider an arbitrary edge $uv$ of $T$ and consider the following subcases.

Let $c(u) = c(v) = 3$. Then the subdivided vertices between $u$ and $v$ must receive the sequence of colors $1,2,1,3,1,2,1,\ldots$. But then the number of subdivided vertices between $u$ and $v$ is odd, a contradiction.

Let $c(u) = c(v) = 2$ and let $w$ be the vertex adjacent to $u$ on the $u,v$-path. Assume first that $c(w) = 1$. Then the vertices between $u$ and $v$ receive colors $1,3,1,2,1,3,1,\ldots$ which again mean that there are an odd number of these subdivided vertices.  Assume next that $c(w) = 3$. We may assume that $y\ne v$ is another neighbor of $u$ in $T$. Then the neighbor of $u$ on the $u,y$-path in $S_i(T)$ receives color $1$. But then we need color at least $4$ for the next vertex on the $u,y$-path.

Suppose finally that $c(u) = 2$ and $c(v) = 3$. If $c(w) = 1$, then the sequence of colors on the $u,v$-path is $1,3,1,2,1,3,\ldots$ and we would have an odd number of subdivided vertices. While if $c(w) = 3$ we get the same contradiction as in the above paragraph.
\qed

\section*{Acknowledgements}

This work was supported in part by the Ministry of Science of Slovenia
under the grant ARRS-BI-US/16-17-013. B.B. and S.K. are also supported in part by the Ministry of Science of Slovenia under the grant P1-0297. D.F.R. is supported by a grant from the Simons Foundation (Grant Number \#209654 to Douglas F. Rall). The authors wish to express their appreciation to Jernej Azarija for the computations of an optimal packing coloring of the subdivided Petersen graph.

\end{document}